\def\dedicace{
\hfill
\vbox{\it
      \hbox{Bon Anniversaire K\`alman:}
      \hbox{tu as 60 ans,}
      \hbox{et on se conna\^{\i}t depuis 30 ans!}
}
}

\magnification=\magstephalf

\def\refBaWu{1}
\def\refCassels{2}
\def\refFeldman{3}
\def\refFN{4}
\def\refLang{5}
\def\refLauMiNe{6}
\def\refMahlerA{7}
\def\refMahlerB{8}
\def\refMatveev{9}
\def\refMignotte{10}
\def\refNW{11}
\def\refPWDurham{12}
\def\refPWSTN{13}
\def\refRW{14}
\def\refWAustralie{15}
\def\refWCanada{16}
\def\refWCrelle{17}
\def\refWDalag{18}
\def\refWielonsky{19}

\def\numeroMahler{(1)}
\def\numeroNW{(2)}
\def\numeroRW{(3)}
\def\numeroThm{(4)}
\def\numeroHauteur{(5)}
\def\numeroLemmeA{(6)}
\def\numeroLemmeB{(7)}
\def\numeroDemThA{(8)}
\def\numeroDemThB{(9)}
\def\numeroDemThC{(10)}
\def\numeroDemThD{(11)}
\def\numeroDemThE{(12)}

\def\date{{\the\day}\
{\ifcase\month\or January\or February\or March\or April
\or May \or June\or July\or August\or September
\or October\or November\or December\fi}\
{\the\year}}
\newcount\hour\newcount\minute	\hour=\time
\divide\hour by 60
\minute=-\hour \multiply\minute by 60 \advance\minute by \time
\def\heure{\the\hour h\kern -.2mm
\ifnum\minute<10 0\fi\the\minute}

\newif\ifpagetitre            \pagetitretrue

\footline{\sevenrm
\ifpagetitre
{\seventt
http{$:$}//www.math.jussieu.fr/${\scriptstyle \sim}$miw/articles/Debrecen.html}
\hfill
\global\pagetitrefalse
\else\ifodd\pageno
\hfill
\else                  \hfill
\fi\fi}

\headline{\sevenrm
\ifpagetitre  To appear in Publ.~Math.~Debrecen,
{\bf 56}, 3-4 (2000).
\hfill\else \hfill\sevenbf \folio\fi}

\font\tenmsb=msbm10
\font\sevenmsb=msbm7
\font\fivemsb=msbm5
\font\seventt=cmtt10 scaled 700
\font\ninett=cmtt10 scaled 900
\font\ninerm=cmr10 scaled 900

\newfam\msbfam
\textfont\msbfam=\tenmsb
\scriptfont\msbfam=\sevenmsb
\scriptscriptfont\msbfam=\fivemsb
\def\Bbb#1{{\fam\msbfam\relax#1}}

\def\cqfd{\unskip\kern 6pt\penalty 500
\raise -2pt\hbox{\vrule\vbox to 10pt{\hrule width 4pt
\vfill\hrule}\vrule}\par}

\def\rmM{{\rm  M}\relax}
\def\rmh{{\rm  h}\relax}

\def\ttA{{\tt A}}
\def\ttB{{\tt B}}
\def\ttI{{\tt I}}
\def\bG{{\Bbb G}}
\def\bQ{{\Bbb Q}}
\def\ttL{{\tt L}}
\def\bR{{\Bbb R}}

\def\bGa{\bG_{\rm a}}
\def\bGm{\bG_{\rm m}}
\def\us{{\underline{s}}}
\def\ut{{\underline{t}}}

\def\uw{\underline{w}}

\def\ubeta{{\underline{\beta}}}
\def\ugamma{{\underline{\gamma}}}
\def\ueta{\underline{\eta}}
\def\cW{{\cal W}}
\def\cV{{\cal V}}
\def\bC{{\Bbb C}}
\def\Card{{\rm  Card }\relax}
\def\cH{{\cal H}}
\def\uT{\underline{T}}
\def\tg#1{T_{#1}(K)}
\def\bP{{\Bbb P}}
\def\bZ{{\Bbb Z}}
\def\ttM{{\tt M}}
\def\LIC{{linear independence condition}}
\def\rank{{\rm  rank }\relax}
\def\hfl#1#2{\smash{\mathop{\hbox to 12 mm{\rightarrowfill}}
\limits^{\scriptstyle#1}_{\scriptstyle#2}}}

\def\virgule{\raise2pt\hbox{,}}

\def\vfl#1#2{\llap{$\scriptstyle #1$}\left\downarrow
\vbox to 6mm{}\right.\rlap{$\scriptstyle #2$}}

\centerline{ \bf On a Problem of Mahler}
\smallskip
\centerline{ \bf  Concerning the
Approximation of Exponentials and Logarithms}
\bigskip

\centerline{\rm by}
\smallskip
\centerline{\it Michel WALDSCHMIDT}

\vskip 1.5 true cm

\dedicace

\vskip 1.5 true cm
\noindent
\underbar{\bf Abstract} We first propose two conjectural
estimates on Diophantine approximation of logarithms of
algebraic numbers. Next we discuss the state of the art and we
give further partial results on this topic.

\vskip 1.5 true cm
\noindent
{\bf \S 1. Two Conjectures on  Diophantine
Approximation of Logarithms of Algebraic Numbers}
\vskip 1 true cm

In 1953 K.~Mahler [\refMahlerA] proved that for any sufficiently
large positive integers $a$ and $b$, the estimates
$$
\Vert \log a\Vert\ge a^{-40\log\log a}
\quad\hbox{and}\quad
\Vert e^b\Vert\ge b^{-40 b}
\eqno\numeroMahler
$$
hold; here, $\Vert\;\cdot\;\Vert$ denotes the distance to the
nearest integer: for $x\in\bR$,
$$
\Vert x\Vert=\min_{n\in\bZ}|x-n|.
$$
In the same paper  [\refMahlerA], he remarks:
\medskip
\item{}{\sl ``The exponent  $40\log\log a$ tends to infinity
very slowly; the theorem is thus not excessively weak, the more
so since one can easily show that
$$
|\log a - b|<{1\over a}
$$
for an infinite increasing sequence of positive integers $a$ and
suitable integers $b$.''}

\medskip\noindent
(We have replaced Mahler's notation $f$ and $a$ by $a$ and $b$
respectively for coherence  with what follows).

\medskip

In view of this remark we shall dub {\it Mahler's problem} the
following open question:

\medskip
\item{($?$)}{\sl Does there exist an absolute constant $c>0$
such that, for any positive integers $a$ and $b$,
$$
|e^b-a|\ge a^{-c}\; ?
$$
}
\medskip\noindent
Mahler's estimates \numeroMahler\ have been refined by Mahler
himself [\refMahlerB], M.~Mignotte [\refMignotte] and
F.~Wielonsky [\refWielonsky]:
 the exponent $40$ can be replaced by $19.183$.

Here we propose two generalizations of Mahler's problem.
One common feature to our two conjectures is that we
replace rational integers by algebraic numbers. However if, for
simplicity, we restrict them to the special case of rational
integers, then they deal with  simultaneous approximation of
logarithms of positive integers by rational integers. In higher
dimension, there are two points of view: one takes either a
hyperplane, or else a line. Our first conjecture is concerned
with lower bounds for $|b_0+b_1\log a_1+\cdots+b_m\log a_m|$,
which amounts to ask for lower bounds for $|e^{b_0}
a_1^{b_1}\cdots a_m^{b_m}-1|$. We are back to the situation
considered by Mahler in the special case $m=1$ and
$b_m=-1$.
Our second conjecture asks for lower bounds for $\max_{1\le i\le
m}|b_i-\log a_i|$, or equivalently for  $\max_{1\le i\le
m}|e^{b_i}- a_i|$. Mahler's problem again corresponds to the
case $m=1$. In both cases $a_1,\ldots,a_m$, $b_0,\ldots,b_m$
are positive rational integers.

Dealing more generally with algebraic numbers, we need to
introduce a notion of height. Here we use Weil's absolute
logarithmic height $\rmh(\alpha)$ (see [\refLang] Chap.~IV,
\S~1, as well as [\refWDalag]), which is related to Mahler's
measure $\rmM(\alpha)$ by
$$
\rmh(\alpha)={1\over d}\log \rmM(\alpha)
$$
and
$$
\rmM(\alpha)=\exp\left(\int_0^1\log|f(e^{2i\pi t})|dt\right),
$$
where $f\in\bZ[X]$ is the minimal polynomial of $\alpha$ and $d$
its degree. Another equivalent definition for
$\rmh(\alpha)$ is given below (\S~3.3).

Before stating our two main conjectures, let us give a special
case, which turns out to be the ``intersection'' of Conjectures
1 and 2 below: it is an extension of Mahler's problem where the
rational integers $a$ and $b$ are replaced by algebraic numbers
$\alpha$ and $\beta$.

\proclaim Conjecture 0. -- There exists a positive absolute
constant  $c_0$ with the following property. Let $\alpha$ and
$\beta$ be complex algebraic numbers and let
$\lambda\in\bC$  satisfy
$e^\lambda=\alpha$. Define $D=[\bQ(\alpha,\beta):\bQ]$.
Further, let $h$ be a positive number satisfying
$$
h\ge  h(\alpha), \quad
h\ge  h(\beta), \quad
h\ge {1\over D} |\lambda|
\quad\hbox{and}\quad
h\ge {1\over D}\cdotp
$$
Then
$$
|\lambda-\beta|\ge
\exp\bigl\{-c_0 D^2h \bigr\}.
$$

One may state this conjecture without introducing the letter
$\lambda$: then the conclusion is  a lower bound for
$|e^\beta-\alpha|$, and the assumption $h\ge |\lambda|/D$ is
replaced by  $h\ge |\beta|/D$. It makes no difference, but for
later purposes we find it more convenient to use logarithms.

\medskip
The best known result in this direction is the following
[\refNW], which includes previous estimates of many authors;
among them are K.~Mahler, N.I.~Fel'dman, P.L.~Cijsouw,
E.~Reyssat, A.I.~Galochkin and G.~Diaz (for references, see
[\refWAustralie], [\refFN], Chap.~2 \S~4.4,
[\refNW]
and
[\refWielonsky]). For convenience we state a
simpler version
\footnote{($^*$)}{\ninerm The main result in  [\refNW] involves
a further parameter $E$ which yields a sharper estimate when
$|\lambda|/D$ is small compared with $h_1$.}

\smallskip
\item{$\bullet$}{\sl Let  $\alpha$ and $\beta$ be  algebraic
numbers and let $\lambda\in\bC$ satisfy
$\alpha=e^\lambda$. Define
$D=[\bQ(\alpha,\beta):\bQ]$. Let $h_1$ and $h_2$ be positive
real numbers satisfying,
$$
h_1\ge h(\alpha),\quad
h_1\ge {1\over D}|\lambda|,\quad
h_1\ge{1\over D}
$$
and
$$
h_2\ge h(\beta),\quad
h_2\ge \log (Dh_1),\quad
h_2\ge \log D,\quad
h_2\ge 1.
$$
Then
$$
|\lambda-\beta|
\ge\exp\Bigl\{-2\cdot 10^{6} D^3h_1h_2(\log D+1) \Bigr\}.
\eqno\numeroNW
$$}

\medskip\noindent
To compare  with Conjecture 0, we notice that from \numeroNW\
we derive, under the assumptions of Conjecture 0,
$$
|\lambda-\beta|\ge
\exp\bigl\{-c D^3h(h+\log D+1)(\log D+1) \bigr\}
$$
with an absolute constant $c$.  This shows how far we are from
Conjecture 0.

In spite of this weakness of the present state of the theory,
we suggest two extensions of Conjecture 0 involving several
logarithms of algebraic numbers. The common hypotheses for our
two conjectures below are the following. We denote by
$\lambda_1,\ldots,\lambda_m$ complex numbers such that the
numbers $\alpha_i=e^{\lambda_i}$ \  ($1\le i\le m$) are
algebraic.  Further, let $\beta_0,\ldots,\beta_m$ be algebraic
numbers. Let $D$  denote  the degree of the number field
$\bQ(\alpha_1,\ldots,\alpha_m,\beta_0,\ldots,\beta_m)$.
Furthermore, let $h$ be a positive number which satisfies
$$
h\ge \max_{1\le i\le m} h(\alpha_i),\quad
h\ge \max_{0\le j\le m} h(\beta_j), \quad
h\ge {1\over D}\max_{1\le i\le m}|\lambda_i|
\quad\hbox{and}\quad
h\ge {1\over D}\cdotp
$$

\proclaim Conjecture 1. -- Assume that the number
$$
\Lambda=\beta_0+\beta_1\lambda_1+\cdots+\beta_m\lambda_m
$$
is non zero.  Then
$$
|\Lambda|\ge
\exp\bigl\{-c_1mD^2h \bigr\},
$$
where $c_1$ is a positive absolute constant.

\proclaim Conjecture 2. --
Assume $\lambda_1,\ldots,\lambda_m$ are linearly independent
over $\bQ$. Then
$$
\sum_{i=1}^m|\lambda_i-\beta_i|\ge
\exp\bigl\{-c_2mD^{1+(1/m)}h \bigr\},
$$
with  a positive absolute constant $c_2$.

\medskip\noindent
{\bf Remark 1.} Thanks to A.O.~Gel'fond,  A.~Baker and others,
a number of results have already been given in the direction of
Conjecture 1. The best known estimates to date are those in
[\refPWDurham], [\refWCanada], [\refBaWu] and [\refMatveev].
Further, in the special case $m=2$, $\beta_0=0$, sharper
numerical values for the constants are known [\refLauMiNe].
However  Conjecture 1 is much stronger than all known lower
bounds:
\smallskip
\item{-} in terms of $h$: best known
estimates involve $h^{m+1}$ in place of $h$;
\smallskip
\item{-} in terms of $D$: so far, we have essentially $D^{m+2}$
in place of $D^2$;
\smallskip
\item{-} in terms of $m$: the sharpest (conditional) estimates,
due to E.M.~Matveev [\refMatveev], display $c^m$ (with an
absolute constant
$c>1$) in place of $m$.
\smallskip\noindent
On the other hand for concrete applications like those
considered by K.~Gy\H{o}ry, a key point is often not to know
sharp estimates in terms of the dependence in the different
parameters, but to have non trivial lower bounds with small
numerical values for the constants. From this point of view a
result like [\refLauMiNe], which deals only with the special
case $m=2$, $\beta_0=0$, plays an important role in many
situations, in spite of the fact that the dependence in the
height of the coefficients $\beta_1, \beta_2$ is not as sharp
as other more general estimates from Gel'fond-Baker's method.

\medskip\noindent
{\bf Remark 2.} In case $D=1$, $\beta_0=0$, sharper estimates
than Conjecture 1 are suggested by Lang-Waldschmidt in
[\refLang], Introduction to Chapters X and XI. Clearly, our
Conjectures 1 and 2 above are not the final word on this topic.

\medskip\noindent
{\bf Remark 3.}
Assume $\lambda_1,\ldots,\lambda_m$ as well as $D$ are fixed
(which means that the absolute constants $c_1$ and $c_2$ are
replaced by  numbers which may depend on $m$,
  $\lambda_1,\ldots,\lambda_m$ and $D$). Then both conjectures
are true: they follow for instance from \numeroNW. The same
holds if
 $\beta_0,\ldots,\beta_m$ and $D$ are fixed.

\medskip\noindent
{\bf Remark 4.} In the special case where
$\lambda_1,\ldots,\lambda_m$   are fixed and
$\beta_0,\ldots,\beta_m$ are restricted to be rational numbers,
Khinchine's Transference Principle (see [\refCassels], Chap.~V)
enables one to relate the two estimates provided by Conjecture
1 and Conjecture 2. It would be interesting to extend and
generalize this transference principle so that one could relate
the two conjectures in more general situations.

\medskip\noindent
{\bf Remark 5.}  The following  estimate has been obtained by
N.I. Feld'man in 1960  (see [\refFeldman],  Th. 7.7 Chap.~7
\S 5); it is the sharpest know result in direction of Conjecture
2 when
$\lambda_1,\ldots,\lambda_m$ are fixed:

\smallskip

\item{$\bullet$}{\sl Under the assumptions of Conjecture 2,
$$
\sum_{i=1}^m|\lambda_i-\beta_i|\ge
\exp\bigl\{-cD^{2+(1/m)}(h+\log D+1)(\log D+1)^{-1}\bigr\}
$$
with  a positive  constant $c$ depending only on
$\lambda_1,\ldots,\lambda_m$.}
\par

\medskip\noindent
Theorem 8.1 in [\refRW] enables one to remove the assumption
that
$\lambda_1,\ldots,\lambda_m$ are fixed, but then yields the
following weaker lower bound:

\item{$\bullet$}{\sl Under the assumptions of Conjecture 2,
$$
\sum_{i=1}^m|\lambda_i-\beta_i|\ge
\exp\bigl\{-cD^{2+(1/m)}h(h+\log D+1)
(\log h+\log D+1)^{1/m}\bigr\},
$$
with  a positive  constant $c$ depending only on $m$.}
\par

\medskip\noindent
As a matter of fact, as in \numeroNW, Theorem 8.1 of [\refRW]
enables one to separate the contribution of the heights of
$\alpha$'s and $\beta$'s.

\item{$\bullet$}{\sl Under the assumptions of Conjecture 2, let
$h_1$ and $h_2$ satisfy
$$
h_1\ge \max_{1\le i\le m} h(\alpha_i), \quad
h_1\ge {1\over D}\max_{1\le i\le m}|\lambda_i|,\quad
h_1\ge {1\over D}
$$
and
$$
h_2\ge \max_{0\le j\le m} h(\beta_j), \quad
h_2\ge \log\log(3Dh_1), \quad
h_2\ge \log D.
$$
Then
$$
\sum_{i=1}^m|\lambda_i-\beta_i|\ge
\exp\bigl\{-cD^{2+(1/m)}h_1 h_2
(\log h_1+\log h_2+2\log D+1)^{1/m}\bigr\},
\eqno\numeroRW
$$
with  a positive  constant $c$ depending only on $m$.}
\par

\medskip\noindent
Again, Theorem 8.1 of [\refRW]  is more precise (it involves the
famous parameter
$E$).

In case $m=1$ the estimate \numeroRW\ gives a lower bound with
$$
D^3h_1 h_2
(\log h_1+\log h_2+2\log D+1),
$$
while \numeroNW\ replaces the factor $(\log h_1+\log h_2+2\log
D+1) $ by $\log D+1$. The explanation of this difference is
that the proof in [\refNW] involves the so-called Fel'dman's
polynomials, while the proof in [\refRW] does not.

\medskip\noindent
{\bf Remark 6.}  A discussion of relations between Conjecture 2
and algebraic independence is given in [\refWDalag], starting
from  [\refRW].

\medskip\noindent
{\bf Remark 7.} One might propose more general conjectures
involving simultaneous linear forms in logarithms. Such
extensions of our conjectures are also suggested by the general
transference principles in [\refCassels]. In this direction a
partial result is given in [\refPWSTN].

\medskip\noindent
{\bf Remark 8.} We deal here with complex algebraic numbers,
which means that we consider only Archimedean absolute values.
The ultrametric situation would be also worth of interest and
deserves to be investigated.

\goodbreak
\vskip 1.5 true cm
\noindent
{\bf \S 2. Simultaneous
Approximation of Logarithms of Algebraic Numbers}
\vskip 1 true cm

Our goal is to give partial results in the direction of
Conjecture 2. Hence we work with several algebraic numbers
$\beta$  (and as many logarithms of algebraic numbers
$\lambda$), but we put them into a matrix $\ttB$. Our
estimates   will be sharper when the rank of $\ttB$ is small.

We need a definition:

\medskip\noindent
{\bf Definition.} A $m\times n$ matrix
$\ttL=(\lambda_{ij})_{1\le i\le m\atop 1\le j\le n}$ satisfies
the {\it \LIC} if, for any non zero tuple
$\ut=(t_1,\ldots,t_m)$ in
$\bZ^m$ and  any non zero tuple $\us=(s_1,\ldots,s_n)$ in
$\bZ^n$, we have
$$
\sum_{i=1}^m\sum_{j=1}^n t_is_j\lambda_{ij}\not=0.
$$

\medskip
This assumption is much stronger than what is actually needed
in the proof, but it is one of the simplest ways of giving a
sufficient condition for our main results to hold.

\proclaim Theorem 1. --
Let $m$, $n$ and $r$ be positive rational integers. Define
$$
\theta={r(m+n)\over mn}\cdotp
$$
There exists a positive constant  $c_1$ with the following
property. Let   $\ttB$ be a
$m\times n$ matrix of rank $\le r$ with coefficients
$\beta_{ij}$ in a number field
$K$.  For $1\le i\le m$ and $1\le j\le n$, let $\lambda_{ij}$ be
a complex number such that the number
$\alpha_{ij}=e^{\lambda_{ij}}$ belongs to $K^\times$ and such
that the $m\times n$ matrix $\ttL=(\lambda_{ij})_{1\le i\le
m\atop 1\le j\le n}$ satisfies the \LIC. Define $D=[K:\bQ]$.
Let $h_1$  and $h_2$ be  positive real numbers satisfying  the
following conditions:
$$
h_1\ge\rmh(\alpha_{ij}) ,\quad
h_1\ge{1\over D}|\lambda_{ij}| , \quad
h_1\ge {1\over D}
$$
and
$$
h_2\ge \rmh(\beta_{ij}),
\quad h_2\ge \log(Dh_1),
\quad h_2\ge \log D, \quad h_2\ge 1
$$
for $1\le i\le m$ and $1\le j\le n$.
Then
$$
\sum_{i=1}^m\sum_{j=1}^n \bigl|\lambda_{ij}-\beta_{ij}\bigr|
\ge
e^{-c_1 \Phi_1}
$$
where
$$
\Phi_1=
\cases{  Dh_1(Dh_2)^\theta & if $Dh_1\ge (Dh_2)^{1-\theta}$,\cr
\mathstrut\cr
(Dh_1)^{1/(1-\theta)} & if $Dh_1<(Dh_2)^{1-\theta}$.
\cr}
\eqno\numeroThm
$$

\medskip\noindent
{\bf Remark 1.}
One could also state the conclusion with the same lower bound
for
$$
\sum_{i=1}^m\sum_{j=1}^n
\bigl|e^{\beta_{ij}}-\alpha_{ij}\bigr|.
$$

\medskip\noindent
{\bf Remark 2.}
Theorem 1 is a variant of Theorem 10.1 in [\refRW]. The main
differences are the following.

In [\refRW], the numbers $\lambda_{ij}$  are fixed (which means
that the final estimate is not explicited in terms of $h_1$).

The second difference is that in [\refRW]
the parameter $r$ is  the rank of the matrix
 $\ttL$. Lemma 1 below shows that
our hypothesis, dealing with the rank of the matrix  $\ttB$,
is less restrictive.

The third difference is that in [\refRW], the linear
independence condition is much weaker than here; but the cost
is that the estimate is slightly weaker in the complex case,
where
$D^{1+\theta}h_2^\theta$ is replaced by
$D^{1+\theta}h_2^{1+\theta}(\log D)^{-1-\theta}$. However it is
pointed out p.~424 of [\refRW] that the conclusion can be
reached  with
$D^{1+\theta}h_2^\theta(\log D)^{-\theta}$
in the special case where all $\lambda_{ij}$ are real number.
It would be interesting to get the sharper estimate without
this extra condition.

Fourthly, the negative power
of $\log D$ which occurs in [\refRW] could be included also in
our estimate by introducing a parameter $E$ (see remark 5
below).

Finally our estimate is sharper than Theorem 10.1 of [\refRW] in case
$Dh_1<(Dh_2)^{1-\theta}$.

\medskip\noindent
{\bf Remark 3.}
In the special case $n=1$, we have $r=1$, $\theta=1+(1/m)$ and
the lower bound
\numeroThm\ is slightly weaker than \numeroRW: according to
\numeroRW, in the estimate
$$
D^{2+(1/m)}h_1h_2^{1+(1/m)},
$$
given by \numeroThm,
one factor $h_2^{1/m}$
can be replaced by
$$
\bigl(\log(eD^2h_1h_2)\bigr)^{1/m}.
$$
Similarly for $n=1$ (by symmetry). Hence Theorem 1 is already
known when
$\min\{m,n\}=1$.

\medskip\noindent
{\bf Remark 4.}
One should stress that \numeroThm\ is
not the sharpest result one can prove.
Firstly the \LIC\ on the matrix $\ttL$ can be weakened.
Secondly the same method enables one to split the dependence of
the different
$\alpha_{ij}$ (see Theorem 14.20 of [\refWDalag]).  Thirdly a
further parameter $E$ can be introduced  (see [\refNW],
[\refWCrelle] and [\refWDalag], Chap.~14 for instance -- our
statement here  corresponds to $E=e$).

\medskip\noindent
{\bf Remark 5.}
In case $Dh_1<(Dh_2)^{1-\theta}$, the number $\Phi_1$ does
not depend on
$h_2$: in fact one does not use the assumption that the numbers
$\beta_{ij}$ are algebraic! Only the rank
$r$ of the matrix comes into the picture.  This follows from
the next result.

\proclaim Theorem 2. --
Let $m$, $n$ and $r$ be positive rational integers with
$mn>r(m+n)$. Define
$$
\kappa={mn\over mn - r(m+n)}\cdotp
$$
There exists a positive constant  $c_2$ with the following
property. Let $\ttL=(\lambda_{ij})_{1\le i\le m\atop
1\le j\le n}$ be a   matrix,
whose entries are logarithms of algebraic numbers, which
satisfies the \LIC.
Let $K$ be a number field containing the algebraic numbers
$\alpha_{ij}=e^{\lambda_{ij}}$ \ ($1\le i\le m$, $1\le j\le
n$). Define $D=[K\colon \bQ]$.  Let $h$
be a positive real number satisfying
$$
h\ge\rmh(\alpha_{ij}),\quad
h\ge{1\over D}|\lambda_{ij}|
\quad\hbox{and}\quad
h\ge {1\over D}
$$
for $1\le i\le m$ and $1\le j\le n$.
Then  for any $m\times n$ matrix  $\ttM=(x_{ij})_{1\le i\le
m\atop 1\le j\le n}$ of rank $\le r$ with complex
coefficients we have
$$
\sum_{i=1}^m\sum_{j=1}^n \bigl|\lambda_{ij}-x_{ij}\bigr|
\ge
e^{-c_2\Phi_2}
$$
where
$$
\Phi_2=(Dh)^\kappa.
$$

Since $\kappa(1-\theta)=1$, Theorem 2 yields the special case of
Theorem 1 where
$Dh_1<(Dh_2)^{1-\theta}$ (cf. Remark 5 above).

\goodbreak
\vskip 1.5 true cm
\noindent
{\bf \S 3. Proofs}
\vskip 1 true cm

Before proving  the theorems, we first deduce \numeroNW\  from
Theorem 4 in [\refNW] and
\numeroRW\ from Theorem 8.1 in [\refRW].

\medskip

The following piece of notation will be convenient: for $n$ and
$S$ positive integers,
$$
\eqalign{
\bZ^n[S]&=[-S,S]^n\cap\bZ^n\cr
&=\bigl\{
\us=(s_1,\ldots,s_n)\in \bZ^n,\; \max_{1\le j\le n}|s_j|\le
S\bigr\}.\cr}
$$
This is a finite set with $(2S+1)^n$ elements.

\goodbreak
\vskip 1  true cm
\noindent
{\bf 3.1.  Proof of \numeroNW}
\vskip .5 true cm

We use Theorem 4 of [\refNW] with $E=e$, $\log A=eh_1$, and we
use the estimates
$$
\rmh(\beta)+\log\max\{1,eh_1\}+\log D+1\le
4h_2\quad\hbox{and}\quad 4e\cdot 105\,500<2\cdot 10^6.
$$
\cqfd

\goodbreak
\vskip 1  true cm
\noindent
{\bf 3.2.  Proof of \numeroRW}
\vskip .5 true cm
We use Theorem 8.1 of [\refRW] with $E=e$, $\log A=eh_1$,
$B'=3D^2h_1h_2$ and $\log B=2h_2$. We may assume without loss
of generality that $h_2$ is sufficiently large with respect to
$m$. The assumption $B\ge D\log B'$ of [\refRW] is satisfied:
indeed the conditions $h_2\ge
\log\log (3Dh_1)$ and $h_2\ge \log D$ imply
$h_2\ge\log\log(3D^2h_1h_2)$.

We need to check
$$
s_1\beta_1+\cdots+s_m\beta_m\not=0 \quad\hbox{for}\quad
\us\in\bZ^m[S]\setminus\{0\}
$$
with
$$
S=\bigl(c_1D\log B')^{1/m}.
$$
Assume on the contrary $s_1\beta_1+\cdots+s_m\beta_m=0$. Then
$$
|s_1\lambda_1+\cdots+s_m\lambda_m|\le mS\max_{1\le i\le m}
|\lambda_i-\beta_i|.
$$
Since $\lambda_1,\ldots,\lambda_m$ are linearly independent, we
may use Liouville's inequality (see for instance [\refWDalag],
Chap.~3) to derive
$$
|s_1\lambda_1+\cdots+s_m\lambda_m|\ge 2^{-D}e^{-mDSh_1}.
$$
In this case one deduces a stronger lower bound than \numeroRW,
with
$$
c D^{2+(1/m)}h_2
\quad\hbox{replaced by}\quad
c' D^{1+(1/m)}.
$$
\cqfd

\goodbreak
\vskip 1  true cm
\noindent
{\bf 3.3.  Auxiliary results}
\vskip .5 true cm

 The proof of the theorems  will require a few preliminary
lemmas.

\proclaim Lemma 1. -- Let  $\ttB=\bigl(\beta_{ij}\bigr)_{1\le
i\le m\atop 1\le j\le n}$ be a matrix whose entries are
algebraic numbers in a field of degree
$D$ and let $\ttL=\bigl(\lambda_{ij}\bigr)_{1\le i\le m\atop
1\le j\le n}$ be a matrix of the same size with complex
coefficients. Assume
$$
\rank (\ttB)>\rank (\ttL).
$$
Let $B\ge 2$ satisfy
$$
\log B\ge \max_{1\le i\le m\atop 1\le j\le n}\rmh(\beta_{ij}).
$$
Then
$$
\max_{1\le i\le m\atop 1\le j\le n}|\lambda_{ij}-\beta_{ij}|
\ge n^{-nD}B^{-n(n+1)D}.
$$

\medskip\noindent
{\it Proof.}
Without loss of generality we may assume that $\ttB$ is a square
regular $n\times n$ matrix. By assumption $\det(\ttL)=0$.

In case $n=1$ we write $\ttB=\bigl(\beta\bigr)$,
$\ttA=\bigl(\lambda\bigr)$ where
$\beta\not=0$ and $\lambda=0$.  Liouville's inequality
([\refWDalag], Chap.~3) yields
$$
|\lambda-\beta|=|\beta|\ge B^{-D}.
$$
Suppose $n\ge 2$. We may assume
$$
\max_{1\le i, j\le n}|\lambda_{ij}-\beta_{ij}|
\le {D\log B\over (n-1)B^D}\virgule
$$
otherwise the conclusion is plain. Since
$$
|\beta_{ij}|\le B^D\quad\hbox{and}\quad
B^{D/(n-1)}\ge 1+{D\over n-1}\log B,
$$
we deduce
$$
\max_{1\le i,j\le n}\max\{|\lambda_{ij}|\; ,\;
|\beta_{ij}|\}\le  B^{nD/(n-1)}.
$$
The polynomial $\det\bigl(X_{ij})$ is homogeneous of degree $n$
and length $n!$; therefore (see Lemma 13.10 of [\refWDalag])
$$
|\Delta|=|\Delta-\det(\ttL)|\le
n\cdot n! \bigl(\max_{1\le i, j\le n}\max\{|\lambda_{ij}|\;
,\;   |\beta_{ij}|\}\bigr)^{n-1} \max_{1\le i, j\le
n}|\lambda_{ij}-\beta_{ij}|.
$$
On the other hand the determinant $\Delta$ of $\ttB$ is a
non zero algebraic number of degree $\le D$. We use Liouville's
inequality again. Now we consider $\det\bigl(X_{ij})$ as a
polynomial of degree $1$ in each of the
$n^2$ variables:
$$
|\Delta|\ge (n!)^{D-1}B^{-n^2D}.
$$
 Finally we conclude the proof of
Lemma 1 by means of the estimate $n\cdot n!\le n^n$.
\cqfd

\medskip

Lemma 1 shows that the assumption $\rank(\ttB)\le r$ of  Theorem
1 is weaker than the condition $\rank(\ttL)=r$ of Theorem 10.1
in [\refRW]. For the proof of   Theorem 1  there is no loss of
generality to assume $\rank(\ttB)= r$ and
$\rank(\ttL)\ge r$.
\bigskip

In the next auxiliary result we use the notion of absolute
logarithmic height on a projective space $\bP_N(K)$, when $K$
is a number field ([\refWDalag], Chap.~3): for
$(\gamma_0:\cdots:\gamma_N)\in \bP_N(K)$,
$$
\rmh(\gamma_0:\cdots:\gamma_N)={1\over D}\sum_{v\in M_K}D_v
\log\max\{|\gamma_0|_v,\ldots,|\gamma_N|_v\},
$$
where $D=[K:\bQ]$, $ M_K$ is the set of normalized absolute
values of $K$, and for
$v\in M_K$, $D_v$ is the local degree.
The normalization of the absolute values is done in such a way
the for $N=1$  we have $\rmh(\alpha)=\rmh(1:\alpha)$.

Here is a simple property of this height.
Let $N$ and $M$ be positive integers and
$\vartheta_1,\ldots,\vartheta_N$, $\theta_1,\ldots,\theta_M$
algebraic numbers. Then
$$
\rmh(1:\vartheta_1:\cdots:\vartheta_N:
\theta_1:\cdots:\theta_M)\le
\rmh(1:\vartheta_1:\cdots:\vartheta_N)+
\rmh(1:\theta_1:\cdots:\theta_M).
$$
One deduces that for algebraic numbers
$\vartheta_0,\ldots,\vartheta_N$, not all of which are zero, we
have
$$
\rmh(\vartheta_0:\cdots:\vartheta_N)\le
\sum_{i=0}^N\rmh(\vartheta_i).
\eqno{\numeroHauteur}
$$

Let $K$ be a number field and  $\ttB$ be a $m\times n$ matrix
of rank $r$ whose entries are in
$K$. There exist two matrices $\ttB'$ and $\ttB''$, of size
$m\times r$ and $r\times n$ respectively, such that
$\ttB=\ttB'\ttB''$. We show how to control the heights of
the entries of $\ttB'$ and $\ttB''$ in terms of the heights of
the entries of $\ttB$  (notice that the proof of Theorem 10.1
in [\refRW] avoids such estimate).

We write
$$
\ttB=\bigl(\beta_{ij}\bigr)_{1\le i\le m\atop 1\le
j\le n},
\quad
\ttB'=\bigl(\beta'_{i\varrho}\bigr)_{1\le i\le m\atop 1\le
\varrho\le r},\quad
\ttB''=\bigl(\beta''_{\varrho j}\bigr)_{1\le \varrho\le r\atop
1\le j\le n}
$$
and we denote by $\ubeta'_1,\ldots,\ubeta'_m$ the $m$ rows of
$\ttB'$ and by
$\ubeta''_1,\ldots,\ubeta''_n$  the $n$ columns of $\ttB''$.
Then
$$
\beta_{ij}=\ubeta'_i\cdot\ubeta''_j \qquad (1\le i\le m, \; 1\le
j\le n),
$$
where  the dot
$\,\cdot\,$ denotes the scalar product in $K^r$.

\proclaim Lemma 2. -- Let $\bigl(\beta_{ij}\bigr)_{1\le i\le
m\atop 1\le j\le n}$ be a $m\times n$ matrix of rank $r$ with
entries in a number field $K$. Define
$$
B=\exp\bigl\{\max_{1\le i\le m\atop 1\le j\le n}\rmh(\beta_{ij}
\bigr\}.
$$
Then there exist
elements
$$
\ubeta'_i=(\beta'_{i1},\ldots,\beta'_{ir}) \quad (1\le i\le m)
\quad\hbox{and}\quad
\ubeta''_j=(\beta''_{1j},\ldots,\beta''_{rj}) \quad (1\le j\le n)
,
$$
in $K^r$ such that
$$
\beta_{ij}=\sum_{\varrho=1}^r
\beta'_{i\varrho}\beta''_{\varrho j}
\quad(1\le i\le m,\; 1\le j\le n)
$$
and such that, for $1\le \varrho\le r$, we have
$$
\rmh(1:\beta'_{1\varrho}:\cdots:\beta'_{m\varrho})\le m\log B
$$
and
$$
\rmh(1:\beta''_{\varrho 1}:\cdots:\beta''_{\varrho n})\le rn\log
B+\log (r!).
\eqno{\numeroLemmeA}
$$

\medskip\noindent
{\it Proof.} We may assume without loss of generality that the
matrix
$\bigl(\beta_{i\varrho}\bigr)_{1\le i,\varrho\le r}$ has rank
$r$. Let $\Delta$ be its determinant. We first take
$\beta'_{i\varrho}=\beta_{i\varrho}$ \ ($1\le i\le m$,
$1\le
\varrho\le r$), so that, by \numeroHauteur,
$$
\rmh(1:\beta'_{1\varrho}:\cdots:\beta'_{m\varrho})\le m\log B
\quad (1\le \varrho\le r).
$$
Next, using Kronecker's symbol, we set
$$
\beta''_{\varrho j}=\delta_{\varrho j} \quad\hbox{for}\quad
  1\le \varrho, j\le r.
$$
Finally we define $\beta''_{\varrho j}$ for $1\le \varrho\le r$,
$r<j\le n$ as the unique solution of the system
$$
\beta_{ij}=\sum_{\varrho=1}^r\beta'_{i\varrho}\beta''_{\varrho
j}
\quad (1\le i\le m,\, r<j\le n).
$$
Then for $1\le \varrho\le r$ we have
$$
(1:\beta''_{\varrho,r+1}:\cdots:\beta''_{\varrho n})=
(\Delta:\Delta_{\varrho,r+1}:\cdots:\Delta_{\varrho n}),
\eqno{\numeroLemmeB}
$$
where, for $1\le \varrho\le r$ and $r<j\le n$, $\Delta_{\varrho
j}$ is (up to sign) the determinant of the $r\times r$ matrix
deduced from the
$r\times (r+1)$ matrix
$$
\pmatrix{
\beta_{11}&\cdots&\beta_{1r}&\beta_{1j}\cr
\vdots    &\ddots&\vdots    &\vdots    \cr
\beta_{r1}&\cdots&\beta_{rr}&\beta_{rj}\cr}
$$
by deleting the $\varrho$-th  column.
From \numeroLemmeB\
one deduces \numeroLemmeA. This completes the proof of Lemma 2.
\cqfd

\bigskip

We need another auxiliary result:

\proclaim Lemma 3. --  Let  $\ttL=\bigl(\lambda_{ij}\bigr)_{1\le
i\le m\atop 1\le j\le n}$ be a  $m\times n$ matrix of complex
numbers which satisfies the
\LIC. Define $\alpha_{ij}=e^{\lambda_{ij}}$ for $i=1,\ldots,m$
and
$j=1,\ldots,n$.
\hfill\break
1) Consider the set
$$
E=\left\{(\ut,\us)\in \bZ^m\times\bZ^n\; ;\; \prod_{i=1}^m
\prod_{j=1}^n
\alpha_{ij}^{t_is_j}=1\right\}.
$$
For each $\us\in\bZ^n\setminus\{0\}$,
$$
\bigl\{\ut\in \bZ^m\; ;\; (\ut,\us)\in E\bigr\}
$$
is a subgroup of $\bZ^m$ of rank $\le 1$, and similarly, for
each
$\ut\in\bZ^m\setminus\{0\}$,
$$
\bigl\{\us\in \bZ^n\; ;\; (\ut,\us)\in E\bigr\}
$$
is a subgroup of $\bZ^n$ of rank $\le 1$.
\hfill\break
2) Fix $\ut\in \bZ^m\setminus\{0\}$. For each positive integer
$S$, the set
$$
\left\{\prod_{i=1}^m \prod_{j=1}^n \alpha_{ij}^{t_is_j}\; ;\;
\us\in\bZ^n[S]\right\}
\subset\bC^\times
$$
has at least $(2S+1)^{n-1}$ elements.

\medskip\noindent
{\it Proof.}
For the proof of 1), fix $\us\in\bZ^n\setminus\{0\}$ and assume
$\ut'$ and $\ut''$ in $\bZ^m$ are such that $(\ut',\us)\in E$
and $(\ut'',\us)\in E$. Taking logarithms we find two rational
integers $k'$ and $k''$ such that
$$
\sum_{i=1}^m \sum_{j=1}^n t'_is_j\lambda_{ij}=2k'\pi\sqrt{-1}
\quad\hbox{and}\quad
\sum_{i=1}^m \sum_{j=1}^n t''_is_j\lambda_{ij}=2k''\pi\sqrt{-1}.
$$
Eliminating $2\pi\sqrt{-1}$ one gets
$$
\sum_{i=1}^m \sum_{j=1}^n(k't''_i-k''t'_i)s_j\lambda_{ij}=0.
$$
Using the \LIC\ on the matrix $\ttL$ one deduces that $\ut'$
and $\ut''$ are linearly dependent over $\bZ$, which proves the
first part of 1). The second part of 1) follows by symmetry.
\smallskip
Now fix   $\ut\in \bZ^m\setminus\{0\}$ and define a mapping
$\psi$ from the finite set $\bZ^n[S]$ to $\bC^\times$ by
$$
\psi(\us)=\prod_{i=1}^m \prod_{j=1}^n \alpha_{ij}^{t_is_j}.
$$
If $\us'$ and $\us''$ in $\bZ^n[S]$ satisfy
$\psi(\us')=\psi(\us'')$, then
$(\us'-\us'',\ut)\in E$. From the first part of the lemma we
deduce that, for each
$\us_0\in
\bZ^n[S]$, the set
$\us-\us_0$, for
$\us$ ranging over the set of elements in $\bZ^n[S]$ for which
$\psi(\us)=\psi(\us_0)$, does not contain two linearly
independent elements. Hence the set
$$
\bigl\{\us\in \bZ^n[S]\; ;\; \psi(\us)=\psi(\us_0)\bigr\}
$$
has at most $ 2S+1$ elements.
Since $\bZ^n[S]$ has $(2S+1)^n$ elements, the conclusion of part
2) of Lemma 3 follows by a simple counting argument (Lemma
7.8    of [\refWDalag]).
\cqfd

\goodbreak
\vskip 1  true cm
\noindent
{\bf 3.4.   Proof of  Theorem 1}
\vskip .5 true cm

As pointed out earlier Theorem 1 in case
$Dh_1<(Dh_2)^{1-\theta}$ is a consequence of Theorem 2 which
will be proved in \S~3.5.  In this section we assume
$Dh_1\ge (Dh_2)^{1-\theta}$ and we prove Theorem 1
with $\Phi_1=Dh_1(Dh_2)^\theta$.

The proof of Theorem 1 is similar to the proof of Theorem 10.1
in [\refRW]. Our main tool is  Theorem 2.1 of [\refWCrelle]. We
do not repeat this statement here, but we check the hypotheses.
For this purpose we need to introduce some notation. We set
$$
d_0=r,\quad
d_1=m,\quad
d_2=0,\quad
d=r+m,
$$
and we consider the algebraic group $G=G_0\times G_1$ with
$G_0=\bGa^r$ and $G_1=\bGm^m$.

There is no loss of generality to assume that the matrix
$\ttB$ has rank
$r$ (since the conclusion is weaker when $r$ is larger).  Hence
we may  use Lemma 2 and introduce the matrix
$$
\ttM=\pmatrix{
           &      &           &\beta''_{11}&\cdots
&\beta''_{1n}\cr
           &\ttI_r&           &\vdots      &\ddots
&\vdots\cr
           &      &           &\beta''_{r1}&\cdots
&\beta''_{rn}\cr
\beta'_{11}&\cdots&\beta'_{1r}&\lambda_{11}&\cdots
&\lambda_{1n}\cr
\vdots     &\ddots&\vdots     &\vdots      &\ddots
&\vdots\cr
\beta'_{m1}&\cdots&\beta'_{mr}&\lambda_{m1}&\cdots
&\lambda_{mn}\cr}
$$
Define $\ell_0=r$ and let $\uw_1,\ldots,\uw_{\ell_0}$ denote the
first $r$ columns of $\ttM$, viewed as elements in $K^{r+m}$:
$$
\uw_k=(\delta_{1k},\ldots,\delta_{rk},
\beta'_{1k},\ldots,\beta'_{mk})
\quad (1\le k\le r)
$$
(with Kronecker's diagonal symbol $\delta$).
The $K$-vector space they span, namely
$W=K\uw_1+\cdots+K\uw_r\subset K^d$, has dimension $r$.

Denote by $\ueta_1,\ldots,\ueta_n$ the last $n$ columns
of $\ttM$, viewed as elements in $\bC^{r+m}$:
$$
\ueta_j=(\beta''_{1j},\ldots,\beta''_{rj},
\lambda_{1j},\ldots,\lambda_{mj})   \quad (1\le j\le n).
$$
Hence for $1\le j\le n$ the point
$$
\ugamma_j=\exp_G\ueta_j=
(\beta''_{1j},\ldots,\beta''_{rj},
\alpha_{1j},\ldots,\alpha_{mj})
$$
lies in $G(K)=K^r\times(K^\times)^m$.

For $\us=(s_1,\ldots,s_n)\in\bZ^n$, define an element
$\ueta_{\us}$ in $\bC^d$ by
$$
\eqalign{
\ueta_{\us}&=s_1\eta_1+\cdots+s_n\eta_n\cr
&=
\left(
\sum_{j=1}^n s_j \beta''_{1j},\ldots,
\sum_{j=1}^n s_j\beta''_{rj},
\sum_{j=1}^n s_j\lambda_{1j},\ldots,
\sum_{j=1}^n s_j\lambda_{mj} \right).
\cr}
$$
Again the point
$$
\gamma_{\us}=\exp_G \ueta_{\us}=
\left(
\sum_{j=1}^n s_j \beta''_{1j},\ldots,
\sum_{j=1}^n s_j\beta''_{rj},
\prod_{j=1}^n \alpha_{1j}^{s_j},\ldots,
\prod_{j=1}^n \alpha_{mj}^{s_j}
\right)
$$
lies in $G(K)$. We denote by
$$
\gamma_{\us}^{(1)}=
\left(
\prod_{j=1}^n \alpha_{1j}^{s_j},\ldots,
\prod_{j=1}^n \alpha_{mj}^{s_j}
\right)\in (K^\times)^m
$$
the projection of $\gamma_{\us}$ on $G_1(K)$.

Next put $\uw'_k=\uw_k$ \ ($1\le k\le r$) and, for  $1\le j\le
n$,
$$
\ueta'_j=(\beta''_{1j},\ldots,\beta''_{rj},\beta_{1j},\ldots,
\beta_{mj})\in K^{r+m},
$$
so that $\uw'_1,\ldots,\uw'_r$, $\ueta'_1,\ldots,\ueta'_n$ are
the column vectors of the matrix
$$
\ttM'=\pmatrix{\ttI_r&\ttB''\cr\ttB'&\ttB\cr}.
$$
Further, for $\us\in\bZ^n$, set
$$
\ueta'_{\us}=s_1\eta'_1+\cdots+s_n\eta'_n.
$$
Consider the vector subspaces
$$
\cW'=\bC \uw'_1+\cdots+\bC \uw'_r
\quad\hbox{and}\quad
\cV'=\bC\ueta'_1+\cdots+\bC\ueta'_n
$$
of  $\bC^d$. Since
$$
\ttM'=
\pmatrix{\ttI_r \cr\ttB'\cr}
\cdot
\Bigl(\ttI_r\; \; \ttB''\Bigr),
$$
the matrix $\ttM'$ has rank $r$, and it follows
that  $\cV'$ and  $\cW'+\cV'$ have dimension $r$.
We  set $r_1=r_2=0$ and
$r_3=r$.

Theorem 2.1 of [\refWCrelle] is completely explicit, it would
not be difficult to derive an explicit value for the constant
$c$ in Theorem 1 in terms of $m$ and $n$ only; but we shall
only show it exists. We denote by $c_0$ a sufficiently large
constant which depend only on
$m$ and
$n$. Without loss of generality we may assume that both $Dh_1$
and
$h_2$ are sufficiently large compared with $c_0$.

We set
$$
S=\left[ (c_0^3Dh_2)^{r/n}\right]
\quad\hbox{and}\quad M=(2S+1)^n,
$$
where the bracket denotes the integral part. Define
$$
\Sigma=\bigl\{\gamma_{\us}\; ;\;
\us\in\bZ^n[S]\bigr\}\subset G(K).
$$
We shall order the elements of $\bZ^n[S]$:
$$
\bZ^n[S]=\bigl\{\us^{(1)},\ldots,\us^{(M)}\bigr\}.
$$
Put $B_1=B_2=e^{c_0h_2}$. The estimates
$$
\rmh\left(1:\sum_{j=1}^n s_j^{(1)} \beta''_{hj}:\cdots:
\sum_{j=1}^n s_j^{(M)}\beta''_{hj}
\right)\le \log B_1
\quad(1\le h\le r)
$$
and
$$
\rmh(1:\beta'_{1k}:\cdots:\beta'_{mk})
\le \log B_2 \quad(1\le k\le r)
$$
follow from Lemma 2 thanks to the conditions  $h_2\ge 1$ and
$h_2\ge\log D$.

Next we set
$$
A_1=\ldots=A_m=\exp\{c_0Sh_1\},\quad E=e.
$$
Thanks to the definition of $h_1$, we have, for $1\le i\le m$,
$$
{e\over D}\le \log A_i,\quad
\rmh\left(\prod_{j=1}^n \alpha_{ij}^{s_j}
\right)\le \log A_i
\quad\hbox{and}\quad
{e\over D}
\left|\sum_{j=1}^n s_j\lambda_{ij}
\right|\le \log A_i.
$$
Then define
$$
T=\left[ (c_0^2Dh_2)^{r/m}\right] ,\quad
V= c_0^{3+4\theta} \Phi_1,\quad
U=V/c_0,\quad
$$
$$
T_0=S_0=  \left[{U\over c_0Dh_2}\right],
\quad
T_1=\cdots=T_m=T,\quad  S_1=\cdots=S_n=S.
$$
The inequalities
$$
DT_0\log B_1\le U,\quad
DS_0\log B_2\le U \quad\hbox{and}\quad
\sum_{i=1}^m DT_i \log A_i\le U
$$
are easy to check.
The integers $T_0,\ldots,T_m$ and $S_0,\ldots,S_n$ are all $\ge
1$, thanks to the assumption
$Dh_1\ge (Dh_2)^{1-\theta}$.
We have $U>c_0D(\log D+1)$ and
$$
{T_0+r\choose r}(T+1)^m>4 V^r.
$$
It will be useful to notice that we also have
$$
S_0^r(2S+1)^n> c_0  T_0^rT^m.
\eqno\numeroDemThA
$$
Finally the inequality
$$
B_2\ge T_0+mT+ dS_0
$$
is satisfied thanks to the conditions $h_2\ge\log(Dh_1)$ and
$h_2\ge\log D$.

Assume now
$$
|\lambda_{ij}-\beta_{ij}|\le e^{-V}
$$
for $1\le i\le m$ and $1\le j\le n$.
Then all hypotheses of Theorem 2.1 of [\refWCrelle] are
satisfied. Hence we obtain an algebraic subgroup
$G^*=G_0^*\times G_1^*$ of $G$, distinct from $G$, such that
$$
S_0^{\ell_0^*}M^* \cH(G^*;\, \uT)\le {(r+m)!\over r!}T_0^rT^m
\eqno{\numeroDemThB}
$$
where
$$
\ell_0^*= \dim_K W^*, \quad
W^*={W+\tg {G^*} \over \tg {G^*}}\virgule \qquad
M^*=\Card (\Sigma^*), \quad
\Sigma^*={\Sigma+G^*(K)\over G^*(K) }\cdotp
$$
Define $d_0^*=\dim(G_0/G_0^*)$ and $d^*=\dim(G/G^*)$.
Since $\cH(G^*;\,\uT)\ge T_0^{r-d_0^*}$, we deduce from
\numeroDemThA\ and
\numeroDemThB
$$
S_0^{\ell_0^*} M^*<S_0^{d_0^*}(2S+1)^n.
\eqno{\numeroDemThC}
$$
We claim $\ell_0^*\ge d_0^*$.  Indeed, consider the diagram
$$
\matrix{
\bC^d    &\hfl{\pi_0}{}  &\bC^r\cr
\vfl{g}{}&               &\vfl{}{g_0}\cr
\bC^{d^*}&\hfl{\pi_0^*}{}&\bC^{d_0^*}\cr}
$$
where
$$
\pi_0:\bC^d\rightarrow \bC^r\quad \hbox{and}\quad
\pi_0^*:\bC^{d^*}\rightarrow \bC^{d_0^*}
$$
denote the
projections with kernels
$$
\{0\}\times \bC^m\quad \hbox{and}\quad
\{0\}\times \tg {G_1^*}
$$
respectively, and
$$
g:\bC^d\rightarrow \bC^{d^*}\quad \hbox{and}\quad
g_0:\bC^r\rightarrow \bC^{d_0^*}
$$
denote the projections
$$
\tg {G}\rightarrow \tg {G}/\tg {G^*}\simeq\tg {G/G^*}
\quad \hbox{and}\quad
\tg {G_0}\rightarrow \tg {G_0}/\tg {G_0^*}\simeq\tg {G_0/G_0^*}
$$
respectively.

We have $W^*=g(W)$ and $\pi_0(W)=\bC^r$. Since $g_0$ is
surjective we deduce
$\pi_0^*(W^*)=\bC^{d_0^*}$, hence
$$
\ell_0^*=\dim W^*\ge \dim\pi_0^*(W^*)=d_0^*.
$$
Combining the inequality $\ell_0^*\ge d_0^*$ with
\numeroDemThC\ we deduce
$$
M^*<(2S+1)^n.
$$
Therefore $\dim G_1^*>0$.
Let $\Sigma_1$
denotes the projection of $\Sigma$ on $G_1$:
$$
\Sigma_1=\left\{
\left(
\prod_{j=1}^n \alpha_{1j}^{s_j},\ldots,
\prod_{j=1}^n \alpha_{mj}^{s_j}
\right)\; ;\;
\us\in\bZ^n[S]\right\}
=\bigl\{\ugamma_{\us^{(1)}}^{(1)},\ldots,
\ugamma_{\us^{(M)}}^{(1)}\bigr\}.
$$
For each $\us'\not=\us''$ in $\bZ^n[S]$ such that
$\ugamma_{\us'}^{(1)}/\ugamma_{\us''}^{(1)}\in
G_1^*(K)$, and for each hyperplane of
$\tg{G^*}$ containing $\tg{G_1^*}$ of equation
$t_1z_1+\cdots+t_mz_m=0$, we get a relation
$$
\prod_{i=1}^m\prod_{j=1}^n\alpha_{ij}^{t_is_j}=1
$$
with $\us=\us'-\us''$.
Using the \LIC\ on  the matrix $\ttL$,
 we
deduce from Lemma 3, part 1), that $G_1^*$ has codimension $1$
in $G_1$; hence
$$
\cH(G^*;\, \uT)
\ge {(r+m-1)!\over r!}
T_0^{r-d_0^*}T^{m-1}.\eqno{\numeroDemThD}
$$
Next from part 2) of Lemma 3 we deduce that the set
$$
\Sigma_1^*={\Sigma_1+G_1^*(K)\over G_1^*(K) }
$$
has at least $(2S+1)^{n-1}$ elements. Hence
$$
M^*=  \Card  (  \Sigma)
\ge  \Card (\Sigma_1^*)
\ge (2S+1)^{n-1} .\eqno{\numeroDemThE}
$$
If $mn\ge m+n$ the estimates  \numeroDemThB, \numeroDemThD\ and
\numeroDemThE\ are not compatible. This contradiction concludes
the proof of Theorem 1 in the case
$\max\{m,n\}>1$ and $Dh_1\ge (Dh_2)^{1-\theta}$. Finally, as we
have seen in Remark 3 of \S~2,   Theorem 1 is already known in
case either $m=1$ or $n=1$.
\cqfd

\goodbreak
\vskip 1  true cm
\noindent
{\bf 3.5.   Proof of  Theorem 2}
\vskip .5 true cm

We start with the easy case where all entries $x_{ij}$ of
$\ttM$ are zero: in this special case Liouville's inequality
 gives
$$
\sum_{i=1}^m\sum_{j=1}^n \bigl|\lambda_{ij}\bigr|
\ge
2^{-D}e^{-Dh}.
$$
Next we remark that we may, without loss of generality,
replace the number $r$ by the actual rank of the matrix $\ttM$.

\par

Thanks to the hypothesis $mn>r(m+n)$, there exist positive
real numbers $\gamma_u$, $\gamma_t$ and $\gamma_s$ satisfying
$$
\gamma_u>\gamma_t+\gamma_s\quad\hbox{and}\quad
r\gamma_u<m\gamma_t<n\gamma_s.
$$
For instance
$$
\gamma_u=1,\quad
\gamma_t={r\over m}+{1\over 2m^2n}\virgule\quad
\gamma_s={r\over n}+{1\over mn^2}
$$
is an admissible choice.

Next let $c_0$ be a sufficiently
large integer. How large it should be can be explicitly
written in terms of $m$, $n$, $r$, $\gamma_u$, $\gamma_t$ and
$\gamma_s$.

 We shall apply Theorem 2.1 of [\refWCrelle] with
$d_0=\ell_0=0$, $d=d_1=m$, $d_2=0$, $G=\bGm^m$,
$r_3=r$, $r_1=r_2=0$,
$$
\ueta_j=(\lambda_{ij})_{1\le i\le m},
\qquad
\ueta'_j=(x_{ij})_{1\le i\le m}\qquad (1\le j\le n).
$$
Since $d_0=\ell_0=0$ we set $T_0=S_0=0$. Therefore the
parameters $B_1$ and $B_2$ will play no role, but for
completenes we set
$$
B_1=B_2=mn(Dh)^{mn}.
$$
We also define $E=e$,
$$
U=c_0^{\gamma_u}(Dh)^\kappa,
\quad V=(12m+9)U,
$$
$$
T_1=\cdots=T_m=T,\quad S_1=\cdots=S_n=S,
$$
where
$$
T=\left[c_0^{\gamma_t}
(Dh)^{r\kappa/m}\right],
\quad
S=\left[ c_0^{\gamma_s}
(Dh)^{r\kappa/n}\right].
$$
Define $A_1=\cdots=A_m$ by
$$
\log A_i={1\over em}c_0^{\gamma_u-\gamma_t-\gamma_s}Sh
\quad(1\le i\le m).
$$
The condition $\gamma_t+\gamma_s<\gamma_u$ enables us to check
$$
\sum_{j=1}^n s_j\rmh(\alpha_{ij})\le \log A_i
\quad\hbox{and}\quad
\sum_{j=1}^n s_j|\lambda_{ij}|\le {D\over E}\log A_i
$$
for $1\le i\le m$ and for any $\us\in \bZ^n[S]$. Moreover,
from the very definition of $\kappa$ we deduce
$$
r\kappa\left({1\over m}+{1\over n}\right)+1=\kappa,
$$
and this  yields
$$
D\sum_{i=1}^m   T_i  \log A_i\le U.
$$
Define
$$
\Sigma=\left\{
\bigl(\alpha_{11}^{s_1}\cdots\alpha_{1n}^{s_n},\ldots,
\alpha_{m1}^{s_1}\cdots\alpha_{mn}^{s_n}\bigr)
\in (K^\times)^m\; ;\;
\us\in \bZ^n[S]\right\}.
$$
From the condition $m\gamma_t>r\gamma_u$ one deduces
$$
(2T+1)^m>2 V^r.
$$
Assume that the conclusion of Theorem 2 does not hold
for $c=c_0^{\gamma_u+1}$.  Then the hypotheses of
Theorem 2.1 of [\refWCrelle]
are satisfied, and we deduce that there exists a connected
algebraic subgroup $G^*$ of $G$, distinct from $G$, which is
incompletely defined by polynomials of multidegrees $\le \uT$
where  $\uT$ stands for the
$m$-tuple $(T,\ldots,T)$,  such that
$$
M^*\cH(G^*;\, \uT)\le m! T^m,
\quad\hbox{ where }\quad
M^*=\Card \left({\Sigma+G^*(K)\over G^*(K)}\right).
$$
Since $m\gamma_t<n\gamma_s$, we
have
$$
m! T^m<(2S+1)^n,
$$
and since  $\cH(G^*;\,  \uT)\ge 1$, we deduce
$$
M^*<(2S+1)^n.
$$
Hence   $ \Sigma[2]\cap
G^*(K)\not=\{e\}$. Therefore there exist
$\us\in\bZ^n[2S]\setminus\{0\}$ and
$\ut\in\bZ^m[T]\setminus\{0\}$ with
$$
\sum_{i=1}^m \sum_{j=1}^n
t_is_j\lambda_{ij}\in 2\pi\sqrt{-1}\bZ.
$$
Let us check, by contradiction, that $G^*$ has codimension
$1$. We already know $G^*\not=G$. If the codimension of
$G^*$ were $\ge 2$, we would have two linearly independent
elements $\ut'$ and $\ut''$ in $\bZ^m[T]$ such that the two
numbers
$$
a'={1\over 2\pi\sqrt{-1}}\sum_{i=1}^m \sum_{j=1}^n
t'_is_j\lambda_{ij}
\quad\hbox{and}\quad
a''={1\over 2\pi\sqrt{-1}}\sum_{i=1}^m \sum_{j=1}^n
t''_is_j\lambda_{ij}
$$
are in $\bZ$. Notice that
$$
\max\{|a'|,|a''|\}\le mnTSDh.
$$
We eliminate $ 2\pi\sqrt{-1}$: set $\ut=a''\ut'-a'\ut''$, so
that
$$
\sum_{i=1}^m \sum_{j=1}^n
t_is_j\lambda_{ij}=0
$$
and
$$
0<|\ut|\le  2mnT^2SDh<(2mnTSDh)^2<U^2.
$$
This is not compatible with our hypothesis that the matrix
$\ttL_{mn}$   satisfies the \LIC.

Hence $G^*$ has codimension $1$ in $G$. Therefore
$$
\cH(G^*;\,  \uT)\ge T^{m-1} \quad\hbox{ and consequently }
\quad
M^*\le m!T.
$$
On the other hand a similar argument shows that any $\us'$,
$\us''$ in
$\bZ^n[2S]$ for which
$$
\sum_{i=1}^m \sum_{j=1}^n
t_is'_j\lambda_{ij}\in 2\pi\sqrt{-1}\bZ
\quad\hbox{and}\quad
\sum_{i=1}^m \sum_{j=1}^n
t_is''_j\lambda_{ij}\in 2\pi\sqrt{-1}\bZ
$$
are linearly dependent over $\bZ$. From Lemma 7.8 of
[\refWDalag] we deduce
$$
M^*\ge  S^{n-1}.
$$
Therefore
$$
S^{n-1}\le m!T .
$$
This is not compatible with the hypotheses $mn>r(m+n)$ and
$r\ge 1$. This final contradiction completes the proof of
Theorem 2.
\cqfd

\goodbreak
\vskip 1.5 true cm
\noindent
{\bf References}
\vskip 1 true cm

\def\refmark#1{\medskip\noindent\hskip -1 true cm
\hbox to .7  true cm{
[#1]\hfill}}

\refmark{\refBaWu}
Baker, Alan; W\"{u}stholz, Gisbert --
Logarithmic forms and group varieties.
J. reine Angew. Math. {\bf 442}  (1993), 19--62.

\refmark{\refCassels}
Cassels, J.W.S. --
{\it An Introduction to Diophantine Approximation.}
Cambridge Tracts in Mathematics and Mathematical Physics, No.
{\bf 45},  Cambridge University Press, New York, 1957.
Reprint of the 1957 edition:   Hafner Publishing Co., New
York, 1972.

\refmark{\refFeldman}
Fel'dman, Naum I. --
{\it Hilbert's seventh problem.}  (Russian)
Moskov. Gos. Univ.,
Moscow, 1982.

\refmark{\refFN}
Fel'dman, Naum I.; Nesterenko, Yuri V. --
{\it Number theory. IV. Transcendental Numbers.}
Encyclopaedia of Mathematical Sciences, {\bf 44}.
Springer-Verlag, Berlin, 1998.

\refmark{\refLang}
Lang, Serge --
{\it Elliptic curves: Diophantine analysis.}
Grundlehren der Mathematischen Wissen\-schaften,
{\bf 231}. Springer-Verlag, Berlin-New York, 1978.

\refmark{\refLauMiNe}
Laurent, Michel; Mignotte, Maurice; Nesterenko, Yuri --
Formes
lin\'{e}aires en deux logarithmes et d\'{e}terminants
d'interpolation.  J. Number Theory {\bf 55} (1995), no.~2,
285--321.

\refmark{\refMahlerA}
Mahler, Kurt --
 On the approximation of logarithms of algebraic
numbers. Philos. Trans. Roy. Soc. London. Ser. A. {\bf 245},
(1953). 371--398.

\refmark{\refMahlerB}
Mahler, Kurt --
Applications of some formulae by Hermite to the
approximation of exponentials and logarithms.
Math. Ann. {\bf 168} (1967) 200--227.

\refmark{\refMatveev}
Matveev, Eug\`ene M. --
Explicit lower estimates for rational homogeneous linear forms
in logarithms of algebraic numbers.
Izv. Akad. Nauk SSSR. Ser. Mat. {\bf 62} No 4, (1998) 81--136.
Engl.~transl.: Izvestiya Mathematics  {\bf 62} No 4, (1998)
723--772.

\refmark{\refMignotte}
Mignotte, Maurice --
Approximations rationnelles de $\pi$ et quelques
autres nombres.  Journ\'ees Arithm\'etiques (Grenoble, 1973),
121--132.  Bull. Soc. Math. France, M\'em. {\bf 37}, Soc. Math.
France, Paris, 1974.

\refmark{\refNW}
Nesterenko, Yuri V.; Waldschmidt, Michel --
On the approximation of the values  of exponential function and
logarithm by algebraic numbers.
(Russian)   {\it Diophantine approximations, Proceedings of
papers dedicated to the memory of Prof. N.~I.~Fel'dman}, ed{.}
Yu.~V.~Nesterenko, Centre for applied research under
Mech.-Math{.} Faculty of MSU, Moscow (1996), 23--42.

\refmark{\refPWDurham}
Philippon, Patrice; Waldschmidt, Michel --
Lower bounds for linear forms
in logarithms.
{\it New advances in transcendence theory} (Durham, 1986),
280--312, Cambridge Univ. Press, Cambridge-New York, 1988.

\refmark{\refPWSTN}
Philippon, Patrice; Waldschmidt, Michel --
Formes lin\'{e}aires de
logarithmes simultan\'{e}es sur les groupes alg\'{e}briques
commutatifs.  {\it S\'{e}minaire de Th\'{e}orie des Nombres,}
Paris 1986--87, 313--347, Progr. Math., {\bf 75},
Birkh\"{a}user Boston, Boston, MA, 1988.

\refmark{\refRW}
Roy, Damien; Waldschmidt, Michel --
Simultaneous approximation and algebraic independence.
The Ramanujan Journal, {\bf 1} Fasc.~4 (1997), 379--430.

\refmark{\refWAustralie}
Waldschmidt, Michel --
Simultaneous  approximation of  numbers  connected  with the
exponential function.  J. Austral. Math. Soc., {\bf25} (1978),
466--478.

\refmark{\refWCanada}
Waldschmidt, Michel --
Minorations de combinaisons lin\'{e}aires de
logarithmes de nombres alg\'{e}\-briques.
Canad. J. Math.  {\bf 45} (1993), no.~1, 176--224.

\refmark{\refWCrelle}
Waldschmidt, Michel --
Approximation diophantienne dans les groupes alg\'{e}briques
commutatifs --- (I) : Une version effective du th\'{e}or\`{e}me
du sous-groupe alg\'{e}brique.
J. reine angew. Math., {\bf 493} (1997), 61--113.

\refmark{\refWDalag}
Waldschmidt, Michel --
{\it Diophantine Approximation on Linear Algebraic Groups.
Transcendence Properties of the Exponential Function in Several
Variables.} Springer Verlag, to appear.
\hfill\break
\hbox{ \seventt
http{$:$}//www.math.jussieu.fr/${\scriptscriptstyle
\sim}$miw/articles/DALAG.html}

\refmark{\refWielonsky}
Wielonsky, Franck --
Hermite-Pad\'e approximants to exponential functions and an
inequality of Mahler. J. Number Theory {\bf 74} (1999), no. 2,
230--249.

\vfill

 \vskip 2truecm plus .5truecm minus .5truecm

\hfill
\vbox{\ninerm
 \hbox{Michel WALDSCHMIDT}
 \hbox{Institut de Math\'ematiques de Jussieu}
 \hbox{Th\'eorie des Nombres\qquad Case 247}
 \hbox{175 rue du Chevaleret}
 \hbox{F--75013 PARIS}
 \hbox{e-mail: {\ninett miw@math.jussieu.fr}}
 \hbox{URL: {\ninett
http{$:$}//www.math.jussieu.fr/${\scriptstyle \sim}$miw/}}
}

\bye